\input amstex
\documentstyle{amsppt}
\define\qtri#1#2#3#4#5{\left[\matrix #1\\ #2 ,#3,#4\endmatrix \right]_{#5}}
\define\qbin#1#2#3{\left[\matrix #1\\ #2\endmatrix \right]_{#3}}
\NoBlackBoxes
\nologo
\topmatter
\title Change of base in Bailey pairs
\endtitle

\author{D. Bressoud,  M. Ismail, and D. Stanton}
\endauthor
\abstract{Versions of Bailey's lemma which change the base from
$q$ to $q^2$ or $q^3$ are given. Iterates of these versions give
many new versions of multisum Rogers-Ramanujan identities. We also prove Melzer's
\cite{7} conjectures for the Fermionic forms of the supersymmetric analogues of
Virasoro characters.}
\endabstract

\endtopmatter

\document

{\subheading{1. Introduction}}

The Bailey chain is a well-known \cite{3} and frequently used technique in the
theory of partitions. It establishes infinite families of equivalent identities,
each identity arising from a Bailey pair which corresponds to a link in the chain.
The Bailey lemma is the recipe for passing between adjacent links. A variation
of the Bailey lemma was described in \cite{1}. It extends the notion of a Bailey chain
to a two-dimensional lattice. The purpose of this paper is to give other explicit
versions of the Bailey lemma which change the base $q$.

This freedom to change the base creates new chains of identities. A wide variety of new
Rogers-Ramanujan identities is the result. For example, iterating the change of base
$q\rightarrow q^2$ yields Theorem 4.3 which, with $a=1$, becomes
$$
\sum_{s_1,\ldots,s_{k+1}} {q^E (-q;q)_{2s_3}(-q^2;q^2)_{2s_4} \cdots
(-q^{2^{k-2}};q^{2^{k-2}})_{2s_{k+1}} \over (q;q)_{s_1-s_2} (q^2;q^2)_{s_2-s_3}
\cdots (q^{2^k};q^{2^k})_{s_{k+1}} }
 =
\prod_{n \not\equiv 0, \pm 2 \ (\text{mod}\ 4+2^k)} \hskip -18pt (1-q^n)^{-1}
$$
where $E = s_1^2 + s_2^2 + s_2 + s_3 + 2s_4 + \cdots + 2^{k-2}s_{k+1}$.

The main theorems are given in \S2. Appropriate limiting cases are stated in
\S3, and
these are used in \S4 to find several new multisum Rogers-Ramanujan
identities.
Using the same techniques in \S5, we verify conjectures of Melzer \cite{7} for the
Fermionic  forms of the supersymmetric analogues of Virasoro characters.
Applications to
basic hypergeometric series are given in \S6. In \S7, we show how to use these
transformations to prove Stembridge's \cite{8} identities of  Rogers-Ramanujan type,
and give a sample of other identities that arise from mixing base changes.

We shall need the definition of a Bailey pair, given below,
and Bailey's lemma, which produces a new Bailey pair from a given Bailey pair.
We use the standard notation found in \cite{5}.

\proclaim{Definition} A pair of sequences $(\alpha_n(a,q),\beta_n(a,q))$
is called a {\it{Bailey pair}} with parameters $(a,q)$ if
$$
\beta_n(a,q)=\sum_{r=0}^n \frac{\alpha_r(a,q)}{(q;q)_{n-r}(aq;q)_{n+r}}
$$
for all $n\ge 0$.
\endproclaim

\proclaim{Bailey's Lemma} Suppose that $(\alpha_n(a,q), \beta_n(a,q))$ is a
Bailey pair
with parameters $(a,q)$. Then $(\alpha_n'(a,q), \beta_n'(a,q))$
is another Bailey pair with
with parameters $(a,q)$, where
$$
\alpha_n'(a,q)=\frac{(\rho_1,\rho_2;q)_n}{(aq/\rho_1,aq/\rho_2;q)_n}
\bigl(\frac{aq}{\rho_1\rho_2}\bigr)^n \alpha_n(a,q),
$$
and
$$
\beta_n'(a,q)=
\sum_{k=0}^n \frac{(\rho_1,\rho_2;q)_k (aq/\rho_1\rho_2;q)_{n-k}}
{(aq/\rho_1,aq/\rho_2;q)_{n}(q;q)_{n-k}}
\bigl(\frac{aq}{\rho_1\rho_2}\bigr)^k\beta_k(a,q).
$$
\endproclaim

{\subheading{2. The main theorems}}

In this section we state and prove versions of Bailey's lemma
in which the base $q$ changes from $q$ to $q^2$ or $q^3$. Theorem 2.2
(Theorem 2.4) is the
inverse of Theorem 2.1 (Theorem 2.3),
and could be considered as changing $q$ to
$q^{1/2}$ ($q^{1/3}$).

\proclaim{Theorem 2.1}
Suppose that $(\alpha_n(a,q), \beta_n(a,q))$ is a Bailey pair
with parameters $(a,q)$. If
$$
\beta_n'(a,q)=
\sum_{k=0}^n
\frac{(-aq;q)_{2k} (B^2;q^2)_k(q^{-k}/B,Bq^{k+1};q)_{n-k}}
{(-aq/B,B;q)_n(q^2;q^2)_{n-k}}
B^{-k}q^{-\binom{k}{2}}\beta_k(a^2,q^2),
$$
then
$(\alpha_n'(a,q), \beta_n'(a,q))$ is a Bailey pair with parameters
$(a,q)$, where
$$
\alpha_r'(a,q)=\frac{(-B;q)_r}{(-aq/B;q)_{r}}B^{-r}
q^{-\binom{r}{2}}\alpha_r(a^2,q^2).
$$
\endproclaim

\demo{Proof} This follows routinely from the definition of a Bailey pair by
interchanging summations and using Singh's quadratic transformation
(III.21) and the $q$-analogue of the Pfaff-Saalsch\"{u}tz theorem (II.12) in
\cite{5}
$$
\ _4\phi_3\left(\left.\matrix q^{-2m},& C^2q^{2m},&D,&Dq\\
&Cq,&Cq^2,&D^2 \endmatrix\right| q^2;q^2\right)=
D^m\frac{(Cq/D,-q;q)_m (1-C)}{(C,-D;q)_m (1-Cq^{2m})}.\tag"(2.1)"
$$
with $m=n-r$, $C=Bq^{-n+2r}$, and $D=-aq^{1+2r}$.
\qed\enddemo
Bailey's lemma is its own inverse, as one could replace
$\rho_1$ and $\rho_2$ by $aq/\rho_1$ and $aq/\rho_2$.
Since Theorem 2.1 changes the base $q$,
its inverse is distinct from Theorem 2.1: Theorem 2.2.

\proclaim{Theorem 2.2}
Suppose that $(\alpha_n(a,q), \beta_n(a,q))$ is a Bailey pair
with parameters $(a,q)$. If
$$
\gamma_n(a,q)=
\sum_{k=0}^n \frac{(qa^2/B;q^2)_{2n-k} (-Bq;q^2)_k}
{(-q^2a^2;q^2)_{2n}(a^4q^2/B^2;q^4)_n(q^4;q^4)_{n-k}}
a^{2k}B^{-k}q^{k^2}\beta_k(a^2,q^2),
$$
then
$(\alpha_n'(a,q), \gamma_n(a,q))$ is a Bailey pair with parameters
$(a^4,q^4)$, where
$$
\alpha_r'(a,q)=\frac{(-Bq;q^2)_r}{(-qa^2/B;q^2)_{r}}
a^{2r}B^{-r} q^{r^2}\alpha_r(a^2,q^2).
$$\endproclaim

\demo{Proof} This follows as in the proof of Theorem 2.1
using the $q$-analogue of the Pfaff-Saalsch\"{u}tz theorem
$$
\align
& _3\phi_2\left(\left.\matrix q^{-2n+2r},& -q^{-2n+2r},&-Bq^{2r+1}\\
&a^2q^{4r+2},&Bq^{1-4n+2r}/a^2 \endmatrix\right|q^2;q^2\right)\\
&\qquad =\ \frac{(-a^2q^{2n+2r},-a^2q^{2r+1}/B;q)_{n-r} }
{(a^2q^{4r+2},a^2q^{2n+1}/B;q)_{n-r} }.
\endalign
$$
\qed\enddemo

\proclaim{Theorem 2.3}
Suppose that $(\alpha_n(a,q), \beta_n(a,q))$ is a Bailey pair
with parameters $(a,q)$. Then $(\alpha_n'(a,q), \beta_n'(a,q))$
is a Bailey pair with parameters $(a^3,q^3)$, where
$$
\alpha_r'(a,q)=a^r q^{r^2}\alpha_r(a,q)
$$
$$
\beta_n'(a,q)=
\frac{1}{(a^3q^3;q^3)_{2n}}\sum_{k=0}^n
\frac{(aq;q)_{3n-k}a^k q^{k^2}}{(q^3;q^3)_{n-k}}
\beta_k(a,q).
\tag"(T1)"
$$
\endproclaim

\demo{Proof} This follows as in the proof of Theorem 2.1
again using  Saalsch\"utz's evaluation \cite{5}
$$
\ _3\phi_2\left(\left.\matrix q^{-n+r},& \omega q^{-n+r},&\omega^2 q^{-n+r}\\
&aq^{2r+1},&q^{r-3n}/a \endmatrix\right|q;q\right)
=\frac{(a\omega q^{r+n+1},a\omega^2 q^{r+n+1};q)_{n-r} }
{(aq^{2r+1},aq^{2n+1};q)_{n-r} }.
$$
where $\omega$ is a primitive cube root of 1.
\qed\enddemo

The inverse of Theorem 2.3 is Theorem 2.4.
\proclaim{Theorem 2.4}
Suppose that $(\alpha_n(a,q), \beta_n(a,q))$ is a Bailey pair
with parameters $(a,q)$. Then $(\alpha_n'(a,q), \beta_n'(a,q))$
is a Bailey pair with parameters $(a,q)$, where
$$
\alpha_r'(a,q)=a^{-r} q^{-r^2}\alpha_r(a^3,q^3)
$$
$$
\align
\beta_n'(a,q)=&
\frac{1}{(aq;q)_{2n}}\sum_{k=0}^n
\frac{(aq^{2n+1};q^{-1})_{3k}(a^3q^3;q^3)_{2(n-k)}}{(q^3;q^3)_{k}}\\
&\times (-1)^{k} q^{3\binom{k}{2}-n^2}a^{-n}\beta_{n-k}(a^3,q^3).
\tag"(T2)"
\endalign
$$
\endproclaim

\demo{Proof} This follows as in the proof of Theorem 2.1
using  the strange $\ _5\phi_4$ evaluation \cite{6, (6.28)}
$$
\sum_{k=0}^m \frac{(q^{-3m};q^3)_k (A^3;q^3)_{2k} q^{3k}}
{(q^3,A^3;q^3)_k (Aq^{1-m};q)_{3k}}=
\frac{q^{-3m^2/2+m/2}(-1)^m (q^3;q^3)_m (1-A)}
{(q^{-1},A^{-1};q^{-1})_m (Aq^{1-m};q)_m (1-Aq^{2m})}
$$
with $A=aq^{1+2r}$ and $m=n-r$.
\qed\enddemo

There is a companion evaluation to (2.1),
which implies a result closely related to Theorem 2.1
$$
\ _4\phi_3\left(\left.\matrix q^{-2m},& C^2q^{2m},&D,&Dq\\
&C,&Cq,&D^2q^2 \endmatrix\right| q^2;q^2\right)=
D^m\frac{(C/D,-q;q)_m}{(C,-Dq;q)_m}.\tag"(2.2)"
$$
We use (2.2) with $m=n-r$, $C=Bq^{-n+2r}$, and $D=-aq^{2r}$ for the next theorem.

\proclaim{Theorem 2.5}
Suppose that $(\alpha_n(a,q), \beta_n(a,q))$ is a Bailey pair
with parameters $(a,q)$. If
$$
\beta_n'(a,q)=
\sum_{k=0}^n
\frac{(-a;q)_{2k} (B^2;q^2)_k(q^{-k+1}/B,Bq^{k};q)_{n-k}}
{(-aq/B,B;q)_n(q^2;q^2)_{n-k}}
B^{-k}q^{k-\binom{k}{2}}\beta_k(a^2,q^2),
$$
then
$(\alpha_n'(a,q), \beta_n'(a,q))$ is a Bailey pair with parameters
$(a,q)$, where
$$
\alpha_r'(a,q)=\frac{(-B;q)_r}{(-aq/B;q)_{r}}\frac{1+a}{1+aq^{2r}}
B^{-r} q^{r-\binom{r}{2}}\alpha_r(a^2,q^2).
$$
\endproclaim

{\subheading{3. Limiting cases}}

It is well-known \cite{3} that Bailey's lemma implies the multisum versions
of the
Rogers-Ramanujan identities due to Andrews.
In this section we record the appropriate limiting cases of Bailey's lemma and
Theorems 2.1-2.4.

First we review \cite{3} the limiting cases of Bailey's lemma which are
used for the Andrews-Gordon identities.
If $\rho_1,\rho_2\rightarrow\infty$ in
Bailey's Lemma, we have
$$
\aligned
\alpha_r'(a,q)&=a^{r} q^{r^2}\alpha_r(a,q),\\
\beta_n'(a,q)&=\ \sum_{k=0}^n \frac{a^k q^{k^2}}{(q;q)_{n-k}}\beta_k(a,q).
\endaligned
\tag S1
$$
Iterate (S1) $k$ times to obtain
$$
\alpha_r^{(k)}(a,q)=a^{rk}q^{kr^2}\alpha_r(a,q). \tag 3.1
$$
If $n\rightarrow\infty$, we have
$$
\align
\beta_{\infty}^{(k)}=&\ \frac{1}{(q;q)_\infty}
\sum_{s_1,\cdots, s_k\ge 0}
\frac{a^{s_1+\cdots+s_k}q^{s_1^2+\cdots+s_k^2}}
{(q;q)_{s_1-s_2}(q;q)_{s_2-s_3}\cdots (q;q)_{s_k}}\beta_{s_k}(a,q)\\
=&\ \frac{1}{(q,aq;q)_\infty}
\sum_{r=0}^\infty a^{rk}q^{kr^2}\alpha_r(a,q).
\tag3.2
\endalign
$$
If we choose the unit Bailey pair \cite{3}
$$
\beta_n(a,q)=
\cases
1, &{\text{if }}n=0\\
0, &{\text{if }}n>0,
\endcases\quad
\alpha_n(a,q)=
\frac{(a;q)_n}{(q;q)_n}\frac{(1-aq^{2n}) }{(1-a)}(-1)^n q^{\binom{n}{2}}
\tag UBP
$$
and then put $a=1$, we obtain a Rogers-Ramanujan identity for modulus $2k+1$
$$
\align
\sum_{s_1,\cdots, s_{k-1}\ge 0}&
\frac{q^{s_1^2+\cdots+s_{k-1}^2}}
{(q;q)_{s_1-s_2}(q;q)_{s_2-s_3}\cdots (q;q)_{s_{k-1}}}\\
=&\ \frac{1}{(q;q)_\infty}\biggl(
1+\sum_{r=1}^\infty q^{(k+1/2)r^2}(q^{-r/2}+q^{r/2})(-1)^r\biggr)\\
=&\ \frac{(q^{2k+1},q^k,q^{k+1};q^{2k+1})_\infty}{(q;q)_\infty}.
\endalign
$$

There are five other choices for iterating Bailey's Lemma which each shift the
modulus of the resulting theta-function by one: If we take $\rho_1\rightarrow\infty,
\rho_2=-\sqrt{aq})$, then we get
$$
\aligned
\alpha_r'(a,q) &=\  a^{r/2} q^{r^2/2}\alpha_r(a,q), \\
\beta_n'(a,q) &=\ \sum_{k=0}^n \frac{(-\sqrt{aq};q)_k}
{(q;q)_{n-k}(-\sqrt{aq};q)_n}a^{k/2} q^{k^2/2}\beta_k(a,q).
\endaligned
\tag S2
$$
When applied to the unit Bailey pair, (S2) has the effect of increasing the
modulus by one instead of 2, in fact (S2)(S2)=(S1).
Thus (S2) may be considered the square root of (S1). If we take $\rho_1\rightarrow\infty,
\rho_2=-q^{1/2}$, then we get
$$
\aligned
\alpha_r'(a,q) &=\ \frac{(-q^{1/2};q)_r}{(-aq^{1/2};q)_r}\, a^{r}
q^{r^2/2}\alpha_r(a,q),
\\
\beta_n'(a,q) &=\ \sum_{k=0}^n \frac{(-q^{1/2};q)_k}
{(q;q)_{n-k}(-aq^{1/2};q)_n}\,a^{k} q^{k^2/2}\beta_k(a,q).
\endaligned
\tag S3
$$
If we take $\rho_1\rightarrow\infty,
\rho_2=-aq^{1/2}$, then we get
$$
\aligned
\alpha_r'(a,q) &=\ \frac{(-a\,q^{1/2};q)_r}{(-q^{1/2};q)_r} 
q^{r^2/2}\alpha_r(a,q),
\\
\beta_n'(a,q) &=\ \sum_{k=0}^n \frac{(-a\,q^{1/2};q)_k}
{(q;q)_{n-k}(-q^{1/2};q)_n}\, q^{k^2/2}\beta_k(a,q).
\endaligned
\tag S4
$$
If we take $\rho_1\rightarrow\infty,
\rho_2=-a^{1/2}q$, then we get
$$
\aligned
\alpha_r'(a,q) &=\ \frac{(-a^{1/2}q;q)_r}{(-a^{1/2};q)_r} 
a^{r/2}q^{(r^2-r)/2}\alpha_r(a,q),
\\
\beta_n'(a,q) &=\ \sum_{k=0}^n \frac{(-a^{1/2}q;q)_k}
{(q;q)_{n-k}(-a^{1/2};q)_n}\, a^{k/2}q^{(k^2-k)/2}\beta_k(a,q).
\endaligned
\tag S5
$$
If we take $\rho_1\rightarrow\infty,
\rho_2=-a^{1/2}$, then we get
$$
\aligned
\alpha_r'(a,q) &=\ \frac{(-a^{1/2};q)_r}{(-a^{1/2}q;q)_r} 
a^{r/2}q^{(r^2+r)/2}\alpha_r(a,q),
\\
\beta_n'(a,q) &=\ \sum_{k=0}^n \frac{(-a^{1/2};q)_k}
{(q;q)_{n-k}(-a^{1/2}q;q)_n}\, a^{k/2}q^{(k^2+k)/2}\beta_k(a,q).
\endaligned
\tag S6
$$
We have that (S1) is the same as (S3)(S4), (S4)(S3), (S5)(S6), or (S6)(S5).

For Theorem 2.1, we have three possible choices of $B$, which change
$\alpha_r(a,q)$
by a quadratic power of $q$,
($B\rightarrow \infty$, $B\rightarrow 0$, and $B^2=aq$).
$$
\aligned
\alpha_r'(a,q)&=\alpha_r(a^2,q^2), \\
\beta_n'(a,q)&=\sum_{k=0}^n \frac{(-aq;q)_{2k}}{(q^2;q^2)_{n-k}}\,
q^{n-k}\beta_k(a^2,q^2),
\endaligned
\tag D1
$$
$$
\aligned
\alpha_r'(a,q)&=a^{-r}q^{-r^2}\alpha_r(a^2,q^2), \\
\beta_n'(a,q)&=\sum_{k=0}^n \frac{(-aq;q)_{2k}}{(q^2;q^2)_{n-k}}\,
q^{k^2+k-2kn-n}(-1)^{n-k}a^{-n}
\beta_k(a^2,q^2),
\endaligned
\tag D2
$$
and
$$
\aligned
\alpha_r'(a,q)&=a^{-r/2}q^{-r^2/2}\alpha_r(a^2,q^2), \\
\beta_n'(a,q)&=\sum_{k=0}^n
\frac{(-aq;q)_{2k}(q^{-1/2-k}/\sqrt{a},q^{k+3/2}\sqrt{a};q)_{n-k}}
{(aq^{2k+1};q^2)_{n-k}(q^2;q^2)_{n-k}}
q^{-\binom{k}{2}}\,(aq)^{-k/2}
\beta_k(a^2,q^2).
\endaligned
\tag D3
$$

For Theorem 2.5 we record only the $B\rightarrow \infty$ case
$$
\aligned
\alpha_r'(a,q)&=\frac{1+a}{1+aq^{2r}}q^r\alpha_r(a^2,q^2), \\
\beta_n'(a,q)&=\sum_{k=0}^n \frac{(-a;q)_{2k}}{(q^2;q^2)_{n-k}}\,
q^{k}\beta_k(a^2,q^2).
\endaligned
\tag D4
$$
The corresponding cases $B \to 0$ and $B^2 = aq$ are labelled (D5) and (D6),
respectively.

Because Theorem 2.3 and Theorem 2.4 have no parameters besides $a$, we label an
application of these
theorems by (T1) and (T2), respectively. We also do not state the analogous
three possibilities for Theorem 2.2 in this paper.

{\subheading{4. Multisum Rogers-Ramanujan identities}}

We have reviewed in \S2 that iterating (S1) gives a multisum
Rogers-Ramanujan identity. In this section and the next section, we consider other
iterates of (S1)--(S6),
(D1)--(D6), (T1)--(T2). We obtain Bressoud's multisum version for even
modulus, and many
new multisum identities.

Before considering the iterates, first we record a proposition
which allows us to insert linear functions of the summation
indices on the multisum side of Rogers-Ramanujan identities.
We need it to change the restricted moduli in (3.2) from
$\equiv 0,\pm k\mod 2k+1$ to $\equiv 0,\pm i \mod 2k+1$. It replaces
the Bailey lattice \cite{1} and is tailored to the choice of $a=1$ in the
unit Bailey pair.

\proclaim{Proposition 4.1}  If $(\alpha_n(q), \beta_n(q))$ is a Bailey pair with
parameters $(1,q)$,
$$
\alpha_n(a,q)=
\cases
1 \text{ for } \qquad n=0,\\
q^{An^2}(q^{(A-1)n}+q^{-(A-1)n})(-1)^n \text{ for }  n>0,
\endcases
$$
then  $(\alpha_n'(q), \beta_n'(q))$ is Bailey pair with parameters $(1,q)$,
where
$\beta_n'(q)=q^n\beta_n(q)$, and
$$
\alpha_n'(q)=
\cases
1 \text{ for } \qquad n=0,\\
q^{An^2}(q^{An}+q^{-An})(-1)^n \text{ for }  n>0.
\endcases
$$
\endproclaim
\demo{Proof} Proposition 4.1 is equivalent to
$$
\sum_{s=-n}^n \qbin{2n}{n-s}{q}w^{s^2-s}(-1)^s=
q^n\sum_{s=-n}^n \qbin{2n}{n-s}{q}w^{s^2-s}(-q)^s
$$
where $w=q^{A}$.
This is easy to verify by considering the $s$ and $1-s$ terms on each side.
\qed\enddemo
We now show how Proposition 4.1 may be used to insert linear factors into
the exponent of $q$ on the sum side of (3.2), thereby changing the excluded
moduli on the product side.
Suppose that we start at the (UBP) with $a=1$,
$$
\alpha_n^{(0)}(q)=q^{n^2/2}(q^{n/2}+q^{-n/2})(-1)^n,
\quad\beta_n^{(0)}(q)=\delta_{0n}.
$$
If we then apply (S1), to obtain a Bailey pair
$(\alpha_n^{(1)}(q),\beta_n^{(1)}(q))$
we have $\alpha_n^{(1)}(q)=q^{3n^2/2}(q^{n/2}+q^{-n/2})(-1)^n$.
We next apply Proposition 4.1 with $A=3/2$ to obtain another Bailey pair
$$\alpha_n^{(2)}(q)=q^{3n^2/2}(q^{3n/2}+q^{-3n/2})(-1)^n,
\quad\beta_n^{(2)}(q)=q^n\beta_n^{(1)}(q).
$$
We could apply (S1) yet again followed by Proposition 4.1 with $A=5/2$,
to obtain
$$\alpha_n^{(4)}(q)=q^{5n^2/2}(q^{5n/2}+q^{-5n/2})(-1)^n,
\quad\beta_n^{(4)}(q)=q^n\beta_n^{(3)}(q).
$$

We see that applying (S1) and Proposition 4.1 alternatively $i$ times inserts
$q^{s_{k-i}+\cdots+s_{k-1}}$ into the left side of (3.2),
and changes the term $q^{-r/2}+q^{r/2}$ on the the right
side to $(q^{-(i+1/2)r}+q^{(i+1/2)r})$. We now have the full form of the
Andrews-Gordon identities,
$$
\align
(q;q)_\infty\beta_{\infty}^{(k)}=&\ 
\sum_{s_1,\cdots, s_{k-1}\ge 0}
\frac{q^{s_1^2+\cdots+s_{k-1}^2+s_{k-i}+\cdots +s_{k-1}}}
{(q;q)_{s_1-s_2}(q;q)_{s_2-s_3}\cdots (q;q)_{s_{k-1}}}\\
=&\ \frac{1}{(q;q)_\infty}
\biggl(1+\sum_{r=1}^\infty
q^{(2k+1)r^2/2}(q^{-(i+1/2)r}+q^{(i+1/2)r})(-1)^r\biggr)\\
=&\ \frac{(q^{2k+1},q^{k-i},q^{k+i+1};q^{2k+1})_\infty}{(q;q)_\infty}.
\endalign
$$

Note that iterating (S1) $k$ times corresponds to adding 2 to the base $k$
times
$$
1@>(S1)>>3@>(S1)>>5@>(S1)>>\cdots@>(S1)>> 2k+1.
$$
For Bressoud's \cite{4} identities of modulus $2k$ we first double the base
using (D1),
then apply (S1) and Proposition 4.1 $i-1$ times, and finally (S1) $k-i$ times,
$$
1@>(D1)>>2@>(S1)>>4@>(S1)>>\cdots@>(S1)>> 2k,
$$

$$
\align
(q;q)_\infty\beta_{\infty}^{(k)}=&\ 
\sum_{s_1,\cdots, s_{k-1}\ge 0}
\frac{q^{s_1^2+\cdots+s_{k-1}^2+s_{k-i}+\cdots +s_{k-1}}}
{(q;q)_{s_1-s_2}(q;q)_{s_2-s_3}\cdots
(q;q)_{s_{k-2}-s_{k-1}}(q^2;q^2)_{s_{k-1}}}\\
=&\ \frac{1}{(q;q)_\infty}
\biggl(1+\sum_{r=1}^\infty q^{kr^2}(q^{-ir}+q^{ir})(-1)^r\biggr)\\
=&\ \frac{(q^{2k},q^{k-i},q^{k+i};q^{2k})_\infty}{(q;q)_\infty}.
\endalign
$$

One may also obtain the modulus $2k$ by using (S1) $k-1$ times and (S2) once
with $a=1$,
$$
1@>(S1)>>3@>(S1)>>5@>(S1)>>\cdots@>(S1)>> 2k-1@>(S2)>> 2k.
$$
By the same method we obtain the generalized G\"ollnitz-Gordon identities
\cite{1, (7.4.4)}
$$
\align
(q;q)_\infty\beta_{\infty}^{(k)}=&\
\sum_{s_1,\cdots, s_{k-1}\ge 0}
\frac{q^{s_1^2/2+\cdots+s_k^2+s_{k-i}+\cdots +s_{k-1}}(-q^{1/2};q)_{s_1}}
{(q;q)_{s_1-s_2}(q;q)_{s_2-s_3}\cdots (q;q)_{s_{k-1}}}\\
=&\ \frac{1}{(q;q)_\infty}
\biggl(1+\sum_{r=1}^\infty q^{kr^2}(q^{-(i+1/2)r}+q^{(i+1/2)r})(-1)^r\biggr)\\
=&\ \frac{(-q;q)_\infty(q^{2k},q^{k-i-1/2},q^{k+i+1/2};q^{2k})_\infty}{(q;q)_\infty}.
\endalign
$$

The Bressoud and G\"ollnitz-Gordon identities may be ``combined"
if we apply (D1) once, (S1) $k-1$ times, and then (S2)
$$
1@>(D1)>>2@>(S1)>>4@>(S1)>>\cdots@>(S1)>> 2k@>(S2)>> 2k+1.
$$
Choosing $a=1$, we obtain
$$
\align
\sum_{s_0,\cdots, s_{k-1}\ge 0}&
\frac{(-q^{1/2};q)_{s_0}q^{s_0^2/2+s_1^2+\cdots+s_{k-1}^2+s_i+\cdots+s_{k-1}}}
{(q;q)_{s_0-s_1}\cdots (q;q)_{s_{k-2}-s_{k-1}}(q^{2};q^{2})_{s_{k-1}}}\\
=&\ \frac{(-q^{1/2};q)_\infty}{(q;q)_\infty}
(q^{2k+1},q^{i+1/2},q^{2k-i+1/2};q^{2k+1})_\infty.
\endalign
$$

Another modulus $2k$ identity may be found by applying (S2)
first and then (S1) $k-1$ times,
$$
1@>(S2)>>2@>(S1)>>4@>(S1)>>\cdots@>(S1)>> 2k
$$
with $a=1$. The result is
$$
\align
(q;q)_\infty\beta_{\infty}^{(k)}=&
\sum_{s_1,\cdots, s_{k-1}\ge 0}
\frac{q^{s_1^2+\cdots+s_{k-1}^2}}
{(q;q)_{s_1-s_2}\cdots (q;q)_{s_{k-2}-s_{k-1}}(q^{1/2};q^{1/2})_{2s_{k-1}}}\\
=&\ \frac{1}{(q;q)_\infty}
\biggl(1+\sum_{r=0}^\infty q^{kr^2}(q^{-r/2}+q^{r/2})\biggr)\\
=&\ \frac{(q^{2k},-q^{k-1/2},q^{k+1/2};q^{2k})_\infty}{(q;q)_\infty}. \tag 4.1
\endalign
$$
This form has an unusual linear perturbation: if we insert
$q^{-(s_i+\cdots+s_{k-1})/2}$, the excluded congruence class does not change,
rather the base changes! Proposition 4.1 does not apply because
only one application of (S2) was used. We state this unusual result in a
proposition.

\proclaim{Proposition 4.2} If $k$ and $i$ are positive integers such that
$1\le i\le k$,
then
$$
\align
\sum_{s_1,\cdots, s_{k-1}\ge 0}&
\frac{q^{s_1^2+\cdots+s_{k-1}^2-(s_i+\cdots+s_{k-1})/2}}
{(q;q)_{s_1-s_2}\cdots (q;q)_{s_{k-2}-s_{k-1}}(q^{1/2};q^{1/2})_{2s_{k-1}}}\\
=&\ \frac{(q^{2i},-q^{i-1/2},-q^{i+1/2};q^{2i})_\infty}{(q;q)_\infty}.
\endalign
$$
\endproclaim
\demo{Proof} This follows immediately from (4.1) and the limiting case of
the $q$-Vandermonde identity (II.7) in \cite{5}
$$
\sum_{s=0}^n \frac{q^{s^2-s/2}}{(q;q)_{n-s}(q^{1/2};q^{1/2})_{2s}}=
\frac{1}{(q^{1/2};q^{1/2})_{2n}}. \quad \qed
$$
\enddemo

If we apply (D1) toward the end, the doubling of the modulus is more
pronounced. For example
$$
1@>(S1)>>3@>(S1)>>\cdots@>(S1)>> 2k-1@>(D1)>>4k-2@>(S1)>>4k
$$
gives for $k\ge 2$, $1\le i\le k$,
$$
\align
\sum_{s_1,\cdots, s_{k}\ge 0}&
\frac{(-q;q)_{2s_2}q^{s_1^2+2s_3^2+\cdots+2s_{k}^2+s_1-s_2+2(s_{i+1}+\cdots+s_k)
}}
{(q^2;q^2)_{s_1-s_2}\cdots (q^2;q^2)_{s_{k-1}-s_{k-2}}(q^{2};q^{2})_{s_{k}}}\\
=&\ \frac{(q^{4k},q^{2i-1},q^{4k-2i+1};q^{4k})_\infty}{(q;q)_\infty}.
\endalign
$$

We may also use (D3) or (D2) instead of (D1). We give two examples using (D3).
The (D3) version of Bressoud's even modulus
theorem is
$$
1@>(D3)>>1@>(S1)>>3@>(S1)>>\cdots@>(S1)>> 2k-1.
$$
If we use $\sqrt{a}=\pm 1$ we have
$$
\align
\sum_{s_1,\cdots, s_{k-1}\ge 0}&
\frac{(\pm q^{-1/2},\pm q^{3/2};q)_{s_{k-1}}q^{s_1^2+\cdots+s_{k-1}^2}}
{(q;q)_{s_1-s_2}\cdots (q;q)_{s_{k-2}-s_{k-1}}(q;q)_{2s_{k-1}}}\\
=&\ \frac{(q^{2k-1},\pm q^{k-3/2},\pm q^{k+1/2};q^{2k-1})_\infty}
{(q;q)_\infty}.
\endalign
$$
A more unusual identity occurs from
$$
1@>(S1)>>3@>(S1)>>\cdots@>(S1)>> 2k-1@>(D3)>>4k-3@>(S1)>>4k-1,
$$
where $k\ge 2$, $1\le i\le k$,
$$
\align
\sum_{s_1,\cdots, s_{k}\ge 0}&
\frac{(-q;q)_{2s_2}(q^{-1/2-s_2},q^{s_2+3/2};q)_{s_1-s_2}
q^{s_1^2-s_2^2/2+2s_3^2+\cdots+2s_{k}^2+2(s_{i+1}+\cdots+s_k)}}
{(q^{2s_2+1};q^2)_{s_1-s_2}(q^{2};q^2)_{s_1-s_2}\cdots
(q^2;q^2)_{s_{k-1}-s_{k-2}}
(q^{2};q^{2})_{s_{k}}}\\
=&\ \frac{(q^{4k-1},q^{2i-3/2},q^{4k-2i+1/2};q^{4k-1})_\infty}{(q;q)_\infty}.
\endalign
$$

Let's take an example
which includes modulus $5$ and uses (D1) $k$ times:
$$
1@>(D1)>>2@>(D1)>>\cdots@>(D1)>> 2^k@>(S1)>> 2^k+2@>(S1)>> 2^k+4.
$$
$$
\alpha_r^{(k+2)}(a,q)=a^{2r}q^{2r^2}\alpha_r(a^{2^k},q^{2^k}).
$$
$$
\align
\beta_n^{(k+2)}(a,q)=&\sum_{s_1,\cdots,s_{k+2}\ge 0}
\frac{a^{s_1+s_2}q^{E}}{(q;q)_{n-s_1}(q;q)_{s_1-s_2}}\\
&\qquad\times \prod_{i=3}^{k+2}
\frac{(-a^{2^{i-3}}q^{2^{i-3}};q^{2^{i-3}})_{2s_i}}
{(q^{2^{i-2}};q^{2^{i-2}})_{s_{i-1}-s_{i}}}
\beta_{s_{k+2}}(a^{2^k},q^{2^k}).
\endalign
$$
where
$$
E=s_1^2+s_2^2+s_2+s_3+2s_4+\cdots+2^{k-2} s_{k+1}-2^{k-1}s_{k+2}.
$$
Choosing the unit Bailey pair (UBP), and letting $n\rightarrow\infty$,
we have the following theorem.
\proclaim{Theorem 4.3} For any non-negative integer $k$,
$$
\align
\frac{1}{(aq;q)_\infty}&\sum_{r=0}^\infty \frac{a^{2r}q^{2r^2+2^k\binom{r}{2}}
(1-a^{2^k}q^{r2^{k+1}})(a^{2^k};q^{2^k})_r(-1)^r}{(1-a^{2^k})(q^{2^k};q^{2^k})_r
}\\
=&\ \sum_{s_1,\cdots,s_{k+1}\ge 0}
\frac{a^{s_1+s_2}q^{E}}{(q;q)_{s_1-s_2}}
\prod_{i=3}^{k+2}
\frac{(-a^{2^{i-3}}q^{2^{i-3}};q^{2^{i-3}})_{2s_i}}
{(q^{2^{i-2}};q^{2^{i-2}})_{s_{i-1}-s_{i}}}
\endalign
$$
where
$$
E=s_1^2+s_2^2+s_2+s_3+\cdots+2^{k-2}s_{k+1} \quad s_{k+2}=0.
$$
\endproclaim
The case $a=1$ of Theorem 4.3 is a Rogers-Ramanujan identity on base $2^k+4$.

\proclaim{Corollary 4.4} For any non-negative integers $k$ and $j$ with
$1\le j \le k$, the generating function for
partitions with part sizes not congruent to $0$ or $\pm(2+2^{j-1})\mod 2^k+4$ is
$$
\sum_{s_1,\cdots,s_{k+1}\ge 0}
\frac{q^{E}}{(q;q)_{s_1-s_2}}
\prod_{i=3}^{k+2}
\frac{(-q^{2^{i-3}};q^{2^{i-3}})_{2s_i}}{(q^{2^{i-2}};q^{2^{i-2}})_{s_{i-1}-s_{i
}}}
=\frac{(q^{2^k+4},q^{2^{j-1}+2},q^{2^k-2^{j-1}+2};q^{2^k+4})_\infty}
{(q;q)_\infty}
$$
where
$$
E=s_1^2+s_2^2+s_2+s_3+\cdots+2^{k-2}s_{k+1}-2^{j-1}s_{j+1}, \quad s_{k+2}=0.
$$
Moreover the same statement holds for part sizes not congruent to $0$ or $\pm 2\mod
2^k+4$, and
not congruent to $0$ or $\pm 1\mod 2^k+4$, if the term $-2^{j-1}s_{j+1}$ in
$E$ is replaced by $0$ or $s_1$, respectively.
\endproclaim

\demo{Proof} The case $E=s_1^2+s_2^2+s_2+s_3+\cdots+2^{k-2}s_{k+1}$ follows
immediately from Theorem 4.3 with $a=1$. We need to insert the appropriate
linear factors via
Proposition 4.1 for the other excluded congruence classes.

To insert $q^{s_1}$, note that after applying (D1) $k$ times and then (S1)
once, we have
$$
\alpha_n^{(k+1)}(1,q)= q^{n^2}q^{2^{k-1}n^2}(q^{2^{k-1}n}+q^{-2^{k-1}n})(-1)^n.
$$
We apply Proposition 4.1 with $A=2^{k-1}+1$ which changes
$\alpha_n^{(k+1)}(1,q)$ to
$$
q^{(2^{k-1}+1)n^2}(q^{(2^{k-1}+1)n}+q^{-(2^{k-1}+1)n})(-1)^n,
$$
then the final application of (S1) gives
$$
\alpha_n^{(k+2)}(1,q)=
q^{(2^{k-1}+2)n^2}(q^{(2^{k-1}+1)n}+q^{-(2^{k-1}+1)n})(-1)^n.
$$
which excludes the classes $0,\pm 1$ by the Jacobi triple product formula.

For the stated values of $j$, we use Proposition 4.1 in reverse to
insert a linear term after $k-j+1$ iterations of (D1). The term
$q^{-2^{j-1}s_{j+1}}$
appears because we use $j-1$ iterations of (D1) after $q^{-s_{j+1}}$ has
been inserted.
\qed\enddemo

Note that for $k=0$ Corollary 4.4 becomes the usual Rogers-Ramanujan
identities for modulus 5.
Thus we have embedded the odd modulus 5 into an infinite family of even
moduli theorems.
Moreover the number of summations for the moduli $2^k+4$ is $k+1$,
compared to  $2^{k-1}+1$ for the known even moduli theorems.

It is natural to ask if there exist other linear perturbations of $E$ in
Corollary 4.4 which will give the missing excluded congruence classes. For
example,
if $k=3$, the classes $0, \pm 5 \mod 12$ do not appear.
However no such perturbation was found for
this case.

{\subheading{5. The Melzer conjectures}}

Melzer \cite{7} conjectured Rogers-Ramanujan multisum representations for some
closely related infinite products. In this section we shall prove the most general
forms of these conjectures
using the methods of \S4. These conjectures arose as generalizations
of expansions of a Fermionic form of the supersymmetric analogues of Virasoro
characters.

\proclaim{Theorem 5.1} For $i = 1,2,\ldots,k$, we have that
$$ \align
& \sum_{n=0}^{\infty} (-1)^na^{kn}q^{kn^2 +
(k-i+1/2)n}\frac{(-q^{1/2};q)_n(-aq^{n+3/2};q)_{\infty}}{(q;q)_n(aq^{n+1};q)_{\infty}}
\\ &
\qquad \times\ \left( 1 + aq^{n+1/2} -
(1+q^{n+1/2})a^iq^{(i-1/2)(2n+1)}\right)
 \\
=&\ \sum_{s_1,\ldots,s_{k-1}=0}^{\infty} \frac{a^{s_1+s_2+\cdots +
s_{k-1}}(-q^{1/2};q)_{s_1}\, q^{s_1^2/2 + s_2^2 + \cdots + s_{k-1}^2 + s_i+\cdots +
s_{k-1}}}{(q;q)_{s_1-s_2} (q;q)_{s_2-s_3} \cdots (q;q)_{s_{k-1}}} \\
=&\ \sum_{s_1,\ldots,s_{2k-2} = 0}^{\infty} \frac{a^{s_1 + s_3 + \cdots +s_{2k-3}}
q^{(s_1^2+s_2^2 + \cdots + s_{2k-2}^2)/2 + s_{2i-1} + s_{2i+1} + \cdots +
s_{2k-3}}}{(q;q)_{s_1-s_2} (q;q)_{s_2-s_3} \cdots (q;q)_{s_{2k-2}}}. 
\endalign $$
\endproclaim

There is a companion identity for which Melzer \cite{7, (2.10)} had only conjectured
the
$a=1$ case.

\proclaim{Theorem 5.2} For $i = 1,2,\ldots,k$, we have that
$$ \align
& \frac{(-a^{1/2}q;q)_{\infty}}{(aq;q)_{\infty}} \sum_{n=0}^{\infty}
(-1)^na^{(k-1/2)n}q^{kn^2 +
(k-i)n}(1-a^iq^{(2n+1)i})\frac{(aq;q)_n}{
(q;q)_n}
 \\
=&\ (-a^{1/2}q;q)_{\infty}\sum_{s_1,\ldots,s_{k-1}=0}^{\infty}
\frac{a^{s_1+s_2+\cdots + s_{k-1}}\, q^{s_1^2 + s_2^2 + \cdots +
s_{k-1}^2 + s_i + s_{i+1}+\cdots + s_{k-1}}}{(q;q)_{s_1-s_2} (q;q)_{s_2-s_3} \cdots
(q;q)_{s_{k-1}}(-a^{1/2}q;q)_{s_{k-1}}}
\\ =&\ \sum_{s_1,\ldots,s_{2k-2} = 0}^{\infty} \frac{a^{(s_1 + s_2 + \cdots
+s_{2k-2})/2} q^{(s_1^2+s_2^2 + \cdots + s_{2k-2}^2)/2 + s_{2i} +
s_{2i+2} + \cdots + s_{2k-2} + S/2}}{(q;q)_{s_1-s_2} (q;q)_{s_2-s_3} \cdots
(q;q)_{s_{2k-2}}}, 
\endalign $$
where $S = s_1-s_2+s_3-s_4 + \cdots + s_{2k-1}$.
\endproclaim

To prove these identities, we need one more Bailey lemma, the one that sits behind
the Bailey lattice and enables us to change the parameter $a$ to $a/q$. It was first
stated and is proven in \cite{1}, lemma 1.2.

\proclaim{Proposition 5.3} Let $(\alpha_n(aq,q),\beta_n(aq,q))$
be a Bailey pair with parameters $(aq,q)$. If
$$
\align
\beta_n'(a,q)=&\ \sum_{k=0}^n \frac{(\rho_1,\rho_2;q)_k (aq/\rho_1\rho_2;q)_{n-k}
(aq/\rho_1\rho_2)^k}{(q;q)_{n-k}(aq/\rho_1,aq/\rho_2;q)_n}
\beta_n(aq,q)\\
\alpha_n'(a,q)=&\ (1-aq)\left(\frac{aq}{\rho_1\rho_2}\right)^n
\frac{(\rho_1,\rho_2;q)_n}{(aq/\rho_1,aq/\rho_2;q)_n} \\ & \qquad \times\ 
\biggl(\frac{\alpha_n(aq,q)}{1-aq^{2n+1}}-aq^{2n-1}\frac{\alpha_{n-1}(aq,q)}{1-aq^{2
n-1}}
\biggr),
\endalign
$$
where $\alpha_{-1}(a,q) = 0$, then $(\alpha_n'(a,q),\beta_n'(a,q))$ is a Bailey
pair with parameters $(a,q)$.
\endproclaim

The proofs of Theorems 5.1 and 5.2 rely on three special cases of this proposition. In
the first, we let $\rho_1$ and $\rho_2$ approach infinity:
$$ \aligned
\beta_n'(a,q)  =&\ \sum_{k=0}^n \frac{a^k q^{k^2}}{(q;q)_{n-k}}\,\beta_k(aq,q), \\
\alpha_n'(a,q) =&\ (1-aq)a^nq^{n^2} \left( \frac{\alpha_n(aq,q)}{1-aq^{2n+1}} -
aq^{2n-1}\frac{\alpha_{n-1}(aq,q)}{1-aq^{2n-1}}\right). 
\endaligned \tag L1 $$
In the second, we let $\rho_1$ approach infinity and set $\rho_2 = -q^{1/2}$:
$$ \aligned
 \beta_n'(a,q) =&\ \sum_{k=0}^n \frac{(-q^{1/2};q)_k a^k
q^{k^2/2}}{(q;q)_{n-k}(-aq^{1/2};q)_n}\,\beta_k(aq,q),
\\
\alpha_n'(a,q) =&\ (1-aq)a^nq^{n^2/2} \frac{(-q^{1/2};q)_n}{(-aq^{1/2};q)_n} \left(
\frac{\alpha_n(aq,q)}{1-aq^{2n+1}} -
aq^{2n-1}\frac{\alpha_{n-1}(aq,q)}{1-aq^{2n-1}}\right). 
\endaligned\tag L2 $$
In the third, we let $\rho_1$ approach infinity and set $\rho_2 = -a^{1/2}q$:
$$ \aligned
 \beta_n'(a,q) =&\ \sum_{k=0}^n \frac{(-a^{1/2}q;q)_k a^{k/2}
q^{(k^2-k)/2}}{(q;q)_{n-k}(-a^{1/2};q)_n}\,\beta_k(aq,q),
\\
\alpha_n'(a,q) =&\ (1-aq)a^{n/2}q^{(n^2-n)/2} \frac{(-a^{1/2}q;q)_n}{(-a^{1/2};q)_n}
\left(
\frac{\alpha_n(aq,q)}{1-aq^{2n+1}} -
aq^{2n-1}\frac{\alpha_{n-1}(aq,q)}{1-aq^{2n-1}}\right). 
\endaligned\tag L3 $$ 
We will use the fact that (L1) is the same as (L2) followed by (S4) which is also (L3)
followed by (S6), a fact that is easily verified by observing their effect on
$\alpha_n(aq,q)$.

To get the multisum in the second line of Theorem 5.1, we start with
$\beta_k^{(0)}(aq,q)$ from the unit Bailey pair (UBP). If $i\geq 3$, then we
apply (S1)
$k-i+1$ times. We then apply (L1) once which changes the parameter $aq$ to $a$, then
apply (S1) $i-3$ times, and finally apply (S3). This yields
$$ \align & \beta_{\infty}^{(k)}(a,q) \\ & = \frac{1}{(q,-aq^{1/2};q)_{\infty}}
 \sum_{s_1,\ldots,s_{k-1}=0}^{\infty} \frac{a^{s_1+\cdots +
s_{k-1}}(-q^{1/2};q)_{s_1}\, q^{s_1^2/2 + \cdots + s_{k-1}^2 + s_i+\cdots +
s_{k-1}}}{(q;q)_{s_1-s_2}  \cdots (q;q)_{s_{k-1}}}.\endalign $$ 
If $i=2$, we apply (S1) $k-1$ times followed by (L2), and if $i=1$ we apply (S1)
$k-1$ times followed by (S3). For purposes of illustration, we assume that $i\geq 3$;
the other cases follow similarly.

We now apply the same sequence of transformations to
$\alpha_n^{(0)}(aq,q)$:
$$ \align
\alpha_n^{(k-i+1)}(aq,q) =&\ (-1)^n a^{(k-i+1)n}q^{(k-i+3/2)n^2 + (k-i+1/2)n}
\frac{(aq;q)_n}{(q;q)_n}\frac{(1-aq^{2n+1})}{(1-aq)}, \\
\alpha_n^{(k)}(a,q) =&\  (1-aq)a^{(i-1)n}q^{(i-3/2)n^2}
\frac{(-q^{1/2};q)_n}{(-aq^{1/2};q)_n}\\ & \qquad \times\ \left(
\frac{\alpha_n^{(k-i+1)}(aq,q)}{1-aq^{2n+1}} -
aq^{2n-1}\frac{\alpha_{n-1}^{(k-i+1)}(aq,q)}{1-aq^{2n-1}}\right).
\endalign $$

It follows that
$$ \align
\beta_{\infty}^{(k)}(a,q) =&\ \frac{1}{(q,aq;q)_{\infty}}
\sum_{n=0}^{\infty}(1-aq)a^{(i-1)n}q^{(i-3/2)n^2}
\frac{(-q^{1/2};q)_n}{(-aq^{1/2};q)_n}\\ & \qquad \times\ \left(
\frac{\alpha_n^{(k-i+1)}(aq,q)}{1-aq^{2n+1}} -
aq^{2n-1}\frac{\alpha_{n-1}^{(k-i+1)}(aq,q)}{1-aq^{2n-1}}\right) \\
=&\ \frac{1}{(q,aq;q)_{\infty}}
\sum_{n=0}^{\infty}\frac{(1-aq)}{(1-aq^{2n+1})}a^{(i-1)n}q^{(i-3/2)n^2}
\alpha_n^{(k-i+1)}(aq,q)
\frac{(-q^{1/2};q)_n}{(-aq^{1/2};q)_n}\\ & \qquad \times\ \left( 1 -
a^iq^{(2n+1)(i-1/2)}\frac{1+q^{n+1/2}}{1+aq^{n+1/2}}
\right) \\
=&\ \frac{1}{(q,aq;q)_{\infty}}
\sum_{n=0}^{\infty}(-1)^n a^{kn}q^{kn^2+(k-i+1/2)n}
\frac{(aq,-q^{1/2};q)_n}{(q;q)_n(-aq^{1/2};q)_{n+1}}
\\ & \qquad \times\ \left( 1+aq^{n+1/2} - a^iq^{(2n+1)(i-1/2)}(1+q^{n+1/2})
\right).
\endalign $$

To get the last multisum of Theorem 5.1, we again start with $\beta_n^{(0)}(aq,q)$
from (UBP), we apply the pair of transformations (S3) followed by (S4) a
total of $k-i+1$ times, then apply (L2), then apply the pair (S4) followed by (S3) a
total of $i-2$ times. If $i=1$, then we just apply (S3)(S4) $k-1$ times, followed
by (S3). Since (S3)(S4) = (S4)(S3) = (S1) and (L2)(S4) = (L1), this is equivalent to
the sequence of transformations used to obtain the first two sums. \qed

A special case of this theorem is Melzer's conjecture (2.6) \cite{7}, a Fermionic
form of the supersymmetric analogue $\hat\chi_{1,2k-2i-1}^{(2,4k)}$, $0\le i\le k-1$
of a Virasoro character:
$$
\align
(q;q)_\infty\beta_{\infty}^{(2k-1)}=&\
\sum_{s_1,\cdots, s_{2k-2}\ge 0}
\frac{q^{s_1^2/2+\cdots+s_{2k-2}^2/2+s_{2k-2i-1}+\cdots+s_{2k-5} +s_{2k-3}}}
{(q;q)_{s_1-s_2}(q;q)_{s_2-s_3}\cdots (q;q)_{s_{2k-2}}}\\
=&\ \frac{(-q^{1/2};q)_\infty}{(q;q)_\infty}
\biggl(1+\sum_{r=1}^\infty q^{kr^2}(q^{-(i+1/2)r}+q^{(i+1/2)r})(-1)^r\biggr)\\
=&\ \frac{(-q^{1/2};q)_\infty(q^{2k},q^{k-i-1/2},q^{k+i+1/2};q^{2k})_\infty}{(q;q)
_\infty}.
\endalign
$$

Theorem 5.2 is proven similarly. We again start with the unit Bailey pair. To get the
summation in the third line, we apply (S2) $2k-2i+1$ times, then (L3), then the
pair (S6)(S5)
$i-2$ times, and finally (S6). This is equivalent to (S2) followed by (S1) $k-i$
times followed by (L1) followed by (S1) $i-2$ times, which can be used to obtain the
summations in the first and second lines.

With $a=1$ in Theorem 5.2, we get the even case of Melzer's (2.6),
$\hat\chi_{1,2i}^{(2,4k)}$:
$$
\align
\sum_{s_1,\cdots, s_{2k-2}\ge 0} &
\frac{q^{s_1^2/2+\cdots+s_{2k-2}^2/2+s_{2i}+\cdots+s_{2k-4} +s_{2k-2}+S/2}}
{(q;q)_{s_1-s_2}(q;q)_{s_2-s_3}\cdots (q;q)_{s_{2k-2}}}\\
=&\ \frac{(-q;q)_\infty}{(q;q)_\infty}
(q^{2k},q^{i},q^{2k-i};q^{2k})_\infty.
\endalign
$$
where $S=(s_1-s_2+s_3-s_4+\cdots +s_{2k-3}-s_{2k-2})$.

Melzer also conjectured \cite{7, (2.3)} alternative forms for
$\hat\chi_{1,2}^{(2,4k)}$ and $\hat\chi_{1,2k}^{(2,4k)}$.
These follow easily in the same way.

{\subheading{6. Basic hypergeometric transformations}}

It is well-known \cite{3} that using Bailey's lemma
twice with the unit Bailey pair gives the terminating version of the
balanced $\ _4\phi_3$ to the very-well poised $\ _8\phi_7$ transformation.
This transformation is a key one in the
theory of basic hypergeometric series. In this section
we record the analogous transformations obtained from Theorem 2.1-2.4 and
Bailey's lemma.
They should be the most important bibasic transformations.

First if we use Bailey's lemma, Theorem 2.1, and the unit Bailey pair we
obtain a
transformation of a balanced $\ _5\phi_4$ to the ``mixed"
very-well poised series
$$
\align
\frac{1}{(aq,q;q)_n}&\sum_{r=0}^n
\frac{(q^{-n},-B;q)_r}{(aq^{n+1},-aq/B;q)_r}
\frac {(\rho_1,\rho_2,a^2;q^2)_r}{
(a^2q^2/\rho_1,a^2q^2/\rho_2,q^2;q^2)_r}
\frac{1-a^2q^{4r}}{1-a^2}\biggl( \frac{a^2q^{n+2}}{B\rho_1\rho_2}\biggr)^r\\
=&\ \frac{(q^{1-n}B;q)_{2n}(-B)^n q^\binom{n}{2}}
{(-aq/B,B;q)_n (q^2;q^2)_n}\\ &\qquad \times\ 
 _5\phi_4\left(\left.\matrix q^{-2n},& B^2,&a^2q^2/\rho_1\rho_2,&-aq,&-aq^2\\
& a^2q^2/\rho_1,&a^2q^2/\rho_2,&Bq^{1-n},&Bq^{2-n} \endmatrix\right|
q^2;q^2\right).
\endalign
$$
This is closely related to \cite{5, (3.10.3)}.

If we first use Theorem 2.1 and Bailey's lemma, and the unit Bailey pair we
obtain
another transformation of a special balanced $\ _5\phi_4$ to the ``mixed"
very-well poised series
$$
\align
&\sum_{r=0}^n
\frac{(q^{-n},-B,\rho_1,\rho_2;q)_r}{(aq^{n+1},-aq/B,aq/\rho_1,aq/\rho_2;q)_r}
\frac{1-a^2q^{4r}}{1-a^2}\biggl(\frac{aq^{n+1}}{B\rho_1\rho_2}\biggr)^r\\
=&\ \frac{(aq,aq/\rho_1\rho_2;q)_{n}}
{(aq/\rho_1,aq/\rho_2;q)_n}
\ _5\phi_4\left(\left.\matrix q^{-n},& Bq,&\rho_1&\rho_2,&1/B&\\
&-aq/B,&B,&\rho_1\rho_2q^{-n}/a,&-q \endmatrix\right| q;q\right).
\endalign
$$

The choice of Bailey's lemma followed by Theorem 2.4 gives yet another
transformation
for a special balanced $\ _6\phi_5$
$$
\align &
\frac{(q^3;q^3)_n (-a)^n q^{3n/2-n^2/2}}{(q;q)_n
(1-aq^{2n+1})(aq^{2-n};q)_{2n-1}} \\ & \qquad \times
\sum_{r=0}^n
\frac{(q^{-n};q)_r}{(aq^{n+1};q)_r}
\frac{(\rho_1,\rho_2,a^3;q^3)_r}{(a^3q^3/\rho_1,a^3q^3/\rho_2,q^3;q^3)_r}
\frac{1-a^3q^{6r}}{1-a^3}\biggl(\frac{a^2q^{n+2}}{\rho_1\rho_2}\biggr)^r\\
&=\ 
 _6\phi_5\left(\left.\matrix
q^{-3n},&
a^3q^3/\rho_1\rho_2,&a^{3/2}q^{3/2}&-a^{3/2}q^{3/2},&a^{3/2}q^{3}&-a^{3/2}q^{3}&
\\
&a^3q^3/\rho_1,&a^3q^3/\rho_2,&aq^{2-n},&aq^{3-n},&aq^{4-n}
\endmatrix\right| q^3;q^3\right).
\endalign
$$

As our final example we take Theorem 2.2 followed by Theorem 2.3 to obtain
$$
\align &
\frac{1}{(a^{12}q^{12};q^{12})_n }
\sum_{r=0}^n
\frac{(q^{-12n};q^{12})_r}{(a^{12}q^{12n+12};q^{12})_r}
\frac{(-Bq,a^2;q^2)_r}{(-qa^2/B,q^2;q^2)_r}
\frac{1-a^2q^{4r}}{1-a^2}\biggl(\frac{a^6q^{12n+5}}{B}\biggr)^r\\
 &=\ 
\frac{(a^4q^4;q^4)_{3n}}{(a^{12}q^{12};q^{12})_{2n}} \\ &\qquad \times
\ _5\phi_4\left(\left.\matrix
q^{-4n},& \omega q^{-4n},&\omega^2 q^{-4n}&qa^2/B,&q^3a^2/B&\\
&q^{-12n}/a^4,&a^4q^2/B^2,&-q^2a^2,&-q^4a^2,& \endmatrix\right| q^4;q^4\right).
\endalign
$$
where $\omega$ is a primitive cube root of 1.

{\subheading{7. Conclusions}}

It is clear that Theorems 2.1--2.4 and Bailey's lemma may be iterated in many
different ways. It is possible to use them to prove all sixteen families of multisum
identities given by Stembridge in \cite{8}. For example, to prove (I14), we start
with the unit Bailey pair with $a=1$, iterate (S1) followed by Proposition 4.1 $k-1$
times, apply (S1) one more time, and then apply (L1) with $a=1/q$ to
get
$$ \align
 (q;q)_{\infty}&
\sum_{s_1,\ldots,s_k \geq 0} {q^{s_1^2 + \cdots + s_k^2 - s_1 + s_2 + \cdots + s_k}
\over (q;q)_{s_1-s_2}
\cdots (q;q)_{s_{k-1}}} \\
=&\ 2 + \sum_{r=1}^{\infty} (q^{(k+3/2)r^2} ( q^r + q^{-r} ) ( q^{(k-1/2)r} +
q^{-(k-1/2)r} ) (-1)^r \\ =&\ (q,q^{2k+2},q^{2k+3};q^{2k+3})_{\infty} +
(q^3,q^{2k},q^{2k+3};q^{2k+3})_{\infty}.\endalign $$  
We subtract
$$ (q;q)_{\infty}
\sum_{s_1,\ldots,s_k \geq 0} {q^{s_1^2 + \cdots + s_k^2 + s_1 + s_2 + \cdots + s_k}
\over (q;q)_{s_1-s_2}
\cdots (q;q)_{s_{k-1}}} = (q,q^{2k+2},q^{2k+3};q^{2k+3})_{\infty}$$
from each side to get Stembridge's (I14):
$$ (q;q)_{\infty}
\sum_{s_1,\ldots,s_k \geq 0} {q^{s_1^2 + \cdots + s_k^2 - s_1 + s_2 + \cdots +
s_k}(1-q^{2s_1}) \over (q;q)_{s_1-s_2}
\cdots (q;q)_{s_{k-1}}} = (q^3,q^{2k},q^{2k+3};q^{2k+3})_{\infty}.$$

For other Rogers-Ramanujan identities, one could consider the
monoid generated by the symbols (S1)--(S6), (D1)--(D6), (T1), and (L1)--(L3)
subject to the relations
$$
\align (S1) =&\ (S2)(S2) = (S3)(S4) = (S4)(S3) = (S5)(S6) = (S6)(S5), \\
(L1) =&\ (L2)(S4) = (L3)(S6), \\
(D1)=&\ (D2)(S2),\\ (D3)=&\ (S2)(D2),\\
(D2)(S2)=&\ (S2)(S2)(D2). \endalign
$$
The number of different representations for a given identity is
the number of words representing a given word.

We state here a few of the identities which may be
obtained from such words. If we take
(D1)(T1)(S1) the result is
$$
\sum_{s_1,s_2\ge 0} \frac{q^{3s_1^2+s_2^2+s_2} (q;q)_{3s_1-s_2}}
{(q^3;q^3)_{2s_1}(q^3;q^3)_{s_1-s_2}(q^2;q^2)_{s_2}}=
\frac{(q^{10},q^4,q^6;q^{10})_\infty}{(q^3;q^3)_\infty}.
\tag7.1
$$

For (S1)(T1)(D1)(S1) we have
$$
\sum_{s_1,s_2\ge 0} \frac{q^{3s_1^2+2s_2^2}
(-q^3;q^6)_{s_1}(q^2;q^2)_{3s_1-s_2}}
{(q^6;q^6)_{2s_1}(q^6;q^6)_{s_1-s_2}(q^2;q^2)_{s_2}}=
\frac{(q^{16},q^7,q^9;q^{16})_\infty}
{(q^{3},q^9,q^{12};q^{12})_\infty}.
\tag7.2
$$

For (S1)(T2)(S1) we have
$$
\sum_{s_1,s_2\ge 0} \frac{q^{3s_1^2+2s_2^2}
(-q^3;q^6)_{s_1}(q^2;q^2)_{3s_1-s_2}}
{(q^6;q^6)_{2s_1}(q^6;q^6)_{s_1-s_2}(q^2;q^2)_{s_2}}=
\frac{(q^{16},q^7,q^9;q^{16})_\infty}
{(q^{3},q^9,q^{12};q^{12})_\infty}.
\tag7.3
$$

For (S1)(T1)(T1)(S1) we have
$$
\sum_{s_1,s_2,s_3\ge 0} \frac{q^{9s_1^2+3s_2^2+s_3^2}
(q^3;q^3)_{3s_1-s_2}(q;q)_{3s_2-s_3}}
{(q^9;q^9)_{2s_1}(q^9;q^9)_{s_1-s_2}(q^3;q^3)_{2s_2}(q^3;q^3)_{s_2-s_3}(q;q)_{s_
3}}=
\frac{(q^{29},q^{14},q^{15};q^{29})_\infty}
{(q^9;q^{9})_\infty}.
\tag7.4
$$

Finally we remark that the new multisums should lead to new
combinatorial interpretations of theorems such as Theorem 4.3.

\address{Mathematics and Computer Science Department,
Macalester College, Saint Paul, Minnesota 55105}
\endaddress
\address{Department of Mathematics, University of South Florida, Tampa,
Florida
33620-5700}
\endaddress
\address{School  of Mathematics, University of Minnesota, Minneapolis,
Minnesota 55455.}
\endaddress

\enddocument
\end